\documentclass[12pt]{article}
\usepackage{graphicx}
\usepackage{amsfonts,mathtools}
\usepackage{epsfig}
\usepackage{multirow}
\usepackage[usenames]{color}
\usepackage{xcolor}
\usepackage{array}
\usepackage{enumerate}
\usepackage{float}

\usepackage{parskip}
\usepackage{amsmath}

\graphicspath{{figures/introduction/eps/}
              {figures/advection/eps/}
              {figures/burgers/eps/}
             }

\newcommand{\be}{\begin{eqnarray}}
\newcommand{\ee}{\end{eqnarray}}
\newcommand{\bes}{\begin{eqnarray*}}
\newcommand{\ees}{\end{eqnarray*}}
\newcommand{\beqn}{\begin{equation}}
\newcommand{\eeqn}{\end{equation}}
\newcommand{\beqns}{\begin{equation*}}
\newcommand{\eeqns}{\end{equation*}}
\newcommand{\bitem}{\begin{itemize}}
\newcommand{\eitem}{\end{itemize}}

\newcommand{\sst}{\scriptscriptstyle}
\newcommand{\ssN}{{\sst N}}

\newcommand{\dt}{{\Delta t}}

\newcommand{\beqa}{\begin{eqnarray}}
\newcommand{\eeqa}{\end{eqnarray}}

\newcommand{\vecu}{{\overrightarrow{u}}}

\newcommand{\bD}  {{\mathbf{D}}}

\newcommand{\half}{\frac 1 2}

\renewcommand{\thetable}{\Roman{table}}

\newcommand{\MajorHead}[1]{\bigskip\begin{center}{\underline{\bf #1}}
                                   \end{center}}

\renewcommand{\thebibliography}[1]{\MajorHead{References}
                  \list {[\arabic{enumi}]}{\settowidth\labelwidth{[#1]}
                         \leftmargin\labelwidth
                         \advance\leftmargin\labelsep
                         \usecounter{enumi}}
                         \def\newblock{\hskip .11em plus .33em minus -.07em}
                         \sloppy \sfcode`\.=1000\relax}

\definecolor{background}{cmyk}{1,1,0,0}

\newtheorem{theorem}{Theorem}

\newtheorem{remark}[theorem]{Remark}

\newcommand{\ssT}{{\sst T}}

\title{Shock Regularization with Smoothness-Increasing Accuracy-Conserving Dirac-Delta Polynomial Kernels}

\author{
B.W. Wissink$^{a,d}$, G.B. Jacobs$^a$, J.K. Ryan$^b$, \\ W.S. Don$^c$ and E.T.A. van der Weide$^d$\\
\textit{$^a$Department of Aerospace Engineering, San Diego State} \\ \textit{University, San Diego, CA}\\
\textit{$^b$School of Mathematics, University of East Anglia,} \\ \textit{Norwich, UK} \\
\textit{$^c$School of Mathematical Sciences, Ocean University} \\ \textit{ of China, Qingdao, China} \\
\textit{$^d$Faculty of Engineering Technology, University of Twente,} \\ \textit{Enschede, The Netherlands} \\
corresponding author, gjacobs@mail.sdsu.edu
}

\begin{document}

\maketitle

\begin{abstract}
A smoothness-increasing accuracy conserving filtering
approach to the regularization of discontinuities is
presented for
single domain spectral collocation approximations of hyperbolic conservation laws.
The filter is based on convolution of a polynomial kernel that approximates
a delta-sequence.  The kernel combines a $k^{th}$ order smoothness with an arbitrary
number  of ${m}$ zero moments. The zero moments ensure a $m^{th}$ order
accurate approximation of the delta-sequence to the delta function.
Through exact quadrature the projection error of the polynomial kernel on the
spectral basis is ensured to be less than the moment error.
A number of test cases on the advection equation, Burger's equation and Euler equations
in 1D and 2D shown that  the filter regularizes
discontinuities while preserving high-order resolution.
\end{abstract}


\noindent
{\bf Keywords} \newline
  Shock capturing, Hyperbolic conservation laws, Regularization, Dirac-Delta, Chebyshev collocation, Filtering
\noindent

\section{Introduction}
\label{chap:INTRODUCTION}

Shock capturing with high-order spectral methods is well known to be plagued by Gibbs phenomena
in the solution. Nonlinear polynomial reconstruction schemes such as the
nonlinear weighted essentially non-oscillatory (WENO) finite difference schemes
on a uniformly spaced grid that have been very successful \cite{Shu} do not
extend well to global polynomial based Chebyshev and Legendre collocation
(spectral) methods. Since spectral methods rely on high-order global basis
functions, slope limiters are often applied to suppress Gibb's oscillations.
The most common approach is to use explicit Runge Kutta
Discontinuous Galerkin (RKDG) methods with min-mod slope limiters introduced by
Cockburn \cite{Cockburn}, \cite{Cockburn2}. A cost-effective alternative to
limiting is the artificial viscosity (AV) approach, e.g. \cite{Chaudhuri},
\cite{Hughes}, \cite{Guermond} where artificial higher even order differential
terms are added to the equations to dissipate the high frequency waves or
smoothen the small scale structures. While this approach is very stable and
more accurate than lower order alternatives, there is no formal proof of
higher-order resolution and accuracy.  Yet another approach is  
to use filtering, which has successfully employed in
simulating shocked flow \cite{don, Chaudhuri}. The high order filter used there is the
variable order exponential filter, which does not satisfy all the criteria for
the definition of filter as laid out in \cite{vandeven} exactly but
asymptotically.

A yet to be explored technique for shock capturing is the use of SIAC-like
filters. The typical application of SIAC filtering is to obtain
superconvergence. This is accomplished by using information that is already
contained in the numerical solution to increase the smoothness of the field and
to reduce the magnitude of the errors. The solution, as a post-processing step,
is convolved against a specifically designed kernel function once at the final
time. The foundations for this postprocessor were established by Bramble and
Schatz \cite{BS}. They showed that the accuracy of Ritz-Galerkin
discretizations can be doubled by convolving the solution against a certain
kernel function. Cockburn et al. \cite{CLSS2} used the ideas of Bramble and
Schatz and those of Mock and Lax \cite{MockLax} to demonstrate that the
postprocessor is also suitable for DG schemes. They proved that, for a certain
class of linear hyperbolic equations with sufficiently smooth solutions, the
postprocessor enhances the accuracy from order $k+1$ to order $2k+1$ in the
$L_2$-norm, where $k$ is the polynomial degree of the original DG
approximation. This postprocessor relies on a symmetric convolution kernel
consisting of $2k+1$ B-splines of order $k+1$.

In \cite{Suarez2} a regularization technique is developed for the
Dirac-$\delta$ source terms in hyperbolic equations that is an excellent
candidate for a kernel of a SIAC-like regularization filter of shock
discontinuities. The technique is based on a class of high-order compactly
supported piecewise polynomials introduced in \cite{Tornberg}. The piecewise
polynomials provide a high-order approximation to the Dirac-Delta whose overall
accuracy is controlled by two conditions: the number of vanishing moments and
smoothness. SIAC kernels have similar smoothness and moment properties as the
Dirac-Delta kernel, but are based on piecewise continuous B-splines instead of
polynomials. In \cite{Suarez3} it was shown that SIAC-like filters based on the
compactly supported Dirac-Delta kernels capture particle-fluid interface
discontinuities with higher-order resolution in single domain spectral
solutions. The compactness of support is closely related to the moment
condition. Smoothness properties yield higher-order resolution away from the
discontinuity.

In the present work we have developed SIAC-like filters based on the high-order
Dirac-Delta kernel for the regularization of shocks and discontinuities. The
filter operation is based on convolution of the solution with the high-order
Dirac-Delta kernels. The formulation of the operation is derived and written in
matrix-vector multiplication form to allow for an efficient and simple
implementation. The filters are tested on a spectral solver for hyperbolic
conservation laws, such as the one dimensional scalar linear advection equation
and scalar nonlinear Burgers' equation, and both one- and two-dimensional
nonlinear Euler equations.  We show that, for the case of the linear advection
equation, filtering discontinuous initial conditions yields a stable,
high-order resolution and with the designed order of accuracy away from the
discontinuity. 

In the area where the solution is the advected (filtered) initial condition,
the designed order of accuracy is achieved and depends on the support width and
number of vanishing moments on the kernel. The support width is made
proportional to grid-spacing and also depends on the number of vanishing
moments.

In the case of the non-linear Burgers' and Euler equations, the filter has to
be applied at every time step in order to maintain stability by smoothing
shocks and discontinuities during the simulation.  This leads to a summation of
filter-errors and smearing of the solution near discontinuities. We show that
by choosing a sufficiently accurate Dirac-Delta kernel and small support-width,
these effects are minimized and design order of accuracy is obtained away
from the discontinuity as well. \newline

In Section \ref{chap:Formulation_and_Methodology} the filter-operation is
derived and background information about SIAC filters, the high-order
Dirac-Delta functions and the Chebyshev collocation method is provided. Next,
Section \ref{chap:Results} presents and discusses the numerical results.
Finally, Section \ref{chap:Conclusions} summarizes the results and gives an
outlook for future work.


\section{Formulation and Methodology} \label{chap:Formulation_and_Methodology}

\subsection{Chebyshev Collocation Method}

Polynomial based spectral methods, in this particular case, the Chebyshev
collocation method, are commonly used in the discretization of spatial
derivatives in PDE's since the order of convergence, for a sufficiently smooth
function, depends only on the smoothness of the function,  also known  as 
spectral accuracy. For example, the spectral approximation of a $C^p$ function
is at least $O( N^{-p})$.  In the case of an analytical function, the
order of convergence of the approximation is exponential \cite{Hesthaven}.  In
the following the Chebyshev collocation method is briefly described for the
purpose of introducing notation. For an overview on spectral methods we refer
to \cite{Hesthaven}.

The collocation method is based on polynomial interpolation of a function $u(x)$, and can be expressed as
\begin{equation}
  u_\ssN(x) = \sum\limits_{j=0}^N u(x_j) l_j(x), \quad l_j(x)=\prod_{k=0,k\neq j}^N \frac{x-x_k}{x_j-x_k} \qquad ,j = 0, \ldots,N,
\label{eq:interpolant}
\end{equation}
where $x_j, j = 0,\ldots,N$ are the collocation points and $l_j(x)$ are the Lagrange interpolation polynomials of degree $N$.
To determine the derivative of the function $u(x)$ at the collocation points $x_i$, $u'(x_i)$, one can simply take the derivative of the Lagrange interpolating polynomial as
\begin{equation}
\frac{\partial u(x_i)}{\partial x} \approx \sum\limits_{j=0}^N u(x_j) l_j'(x_i).
\end{equation}
or, written compactly in the matrix-vector multiplication form as
\begin{equation}
{\vecu}' = \bD \vecu,
\label{eq:spatialdisc}
\end{equation}
where the differentiation matrix $D_{i,j}=l_j'(x_i)$.
In the case of the Chebyshev collocation method, the collocation points are
\begin{equation}
  x_i=-\cos(i\pi/N), \quad i=0,\ldots,N .
\end{equation}
To integrate the resulting system of ordinary differential equations (ODE) in time, we employ the third order Total Variation Diminishing (TVD) Runge-Kutta scheme
\begin{eqnarray}
  u^{(1)}&=&u^n+\Delta t L(u^n)  \nonumber \\
  u^{(2)}&=&\frac{3}{4}u^n+\frac{1}{4}u^{(1)}+\frac{1}{4}\Delta t L(u^{(1)}) \\
  u^{n+1}&=&\frac{1}{3}u^n+\frac{2}{3}u^{(2)}+\frac{2}{3}\Delta t L(u^{(2)}). \nonumber
\end{eqnarray}

\subsection{Dirac-Delta Approximation} \label{chap:Delta-Dirac Approximation}

\newcommand{\veps}{\varepsilon}
\newcommand{\mk}{{m,k}}
\newcommand{\dmk} {\delta_{\veps}^\mk}
\newcommand{\dmkx}{\dmk(x)}
\newcommand{\Pmkx}{P^\mk(x)}
\newcommand{\xeps}{{\Omega^\veps}}
\newcommand{\xepsi}{{\Omega^\veps_i}}
\newcommand{\xepsj}{{\Omega^\veps_j}}
\newcommand{\bS}{\mathbf{S}}

In \cite{Suarez2}, a regularization technique based on a class of high-order,
compactly supported piecewise polynomials is developed that regularizes the
time-dependent, singular Dirac-Delta sources in spectral approximations of
hyperbolic conservation laws. This regularization technique provides
higher-order accuracy away from the singularity.

For the purpose of clarity in the following discussion, we shall define the
compact support domain as $\xeps = [\veps-x, x+\veps]$ and $\xepsi =
[x_i-\veps, x_i+\veps]$ centered at $x=x_i$ with a given support width $\veps >
0$.

\begin{figure}[htbp]
\begin{center}
  \mbox{
        \includegraphics[width=0.40\linewidth]{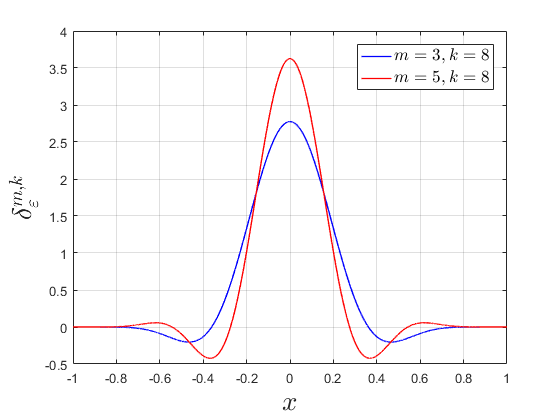}
       }
\end{center}
\caption{Dirac-Delta kernel $\delta_{\varepsilon}^{m,k}$ for $(m,k) = (3,8)$ and $(m,k)=(5,8)$ with scaling parameter $\varepsilon=1$.}
\label{fig:kernel1D}
\end{figure}

We shall then define the Dirac-Delta polynomial kernel, denoted as
$\delta_{\varepsilon}^{m,k}(x)$, that is an approximation of the Dirac-Delta
function $\delta(x)$,

\begin{equation}
\delta_{\varepsilon}^{m,k}(x) = \begin{cases}
    \frac{1}{\varepsilon}P^{m,k}\left(\frac{x}{\varepsilon}\right)       & x \in \xeps, \\
    0  &  \mbox{else},
		\label{eq:delta_function}
  \end{cases}
\end{equation}
where $\varepsilon > 0$ is the scaling parameter.
The function is generated by a multiplication of two lower degree polynomials
that control the number of vanishing moments $m$, and the number of continuous
derivatives at the endpoints of the compact support $k$, respectively.  The
$M=m+2(k+1)$ degree polynomial $P^{m,k}(x)$ is uniquely determined by the
following properties
\begin{eqnarray}
\left( P^{m,k} \right)^{(i)} (\pm 1) &=& 0, \quad i=1,...,k, \label{a1} \\
\int_{-1}^1 P^{m,k}(\xi) d\xi        &=& 1,                  \label{a2} \\
\int_{-1}^1 \xi^i P^{m,k}(\xi) d\xi  &=& 0, \quad i=1,...,m, \label{a3}
\end{eqnarray}
in which (\ref{a1}) determines the number of continuous derivatives at the endpoints ($k$),
(\ref{a2}) states that the function integrates to unity as a Dirac-Delta function, and
(\ref{a3}) determines the number of vanishing moments ($m$).
In Fig. \ref{fig:kernel1D}, we show the polynomial approximation of the Dirac-Delta kernel $\delta_{\varepsilon}^{m,k}(x)$ with $(m,k) = (3,8)$ and $(m,k) = (5,8)$ with scaling parameter $\varepsilon=1$.
It has been shown that the Dirac-Delta approximation $\delta_{\varepsilon}^{m,k}(x)$ has an accuracy of $O(\varepsilon^{m+1})$.

The procedure for the generation of the polynomials $P^{m,k}(\xi)$ is described
in \cite{Tornberg}. A few examples for low values of $m$ and $k$ are given
below.

\begin{itemize}
\item The polynomials with one vanishing moment $m=1$ and with $k=0,1$ and $2$ continuous derivatives, respectively, are
\begin{align}
P^{1,0}=\frac{3}{4}(1-\xi^2), \quad P^{1,1}=\frac{15}{16}(1-2\xi^2+\xi^4), \quad P^{1,2}=\frac{35}{32}(1-3\xi^2+3\xi^4-\xi^6).
\end{align}
\item The polynomials with three vanishing moments $m=3$ and with $k=0,1$ and 2 continuous derivatives, respectively, are
\begin{align}
P^{3,0}&=\frac{15}{32}(3-10\xi^2+7\xi^4), \quad P^{3,1}=\frac{105}{64}(1-5\xi^2+7\xi^4-3\xi^6), \\P^{3,2}&=\frac{315}{512}(3-20\xi^2+42\xi^4-36\xi^6+11\xi^8).
\end{align}
\end{itemize}

\subsection{SIAC Filtering}

The Smoothness-Increasing Accuracy-Conserving (SIAC) filter was developed in the context of superconvergence extraction and error reduction.
SIAC filtering \cite{RSA,kingmirzaee,CKRS,JIRyan} exploits the idea of
superconvergence in the underlying method to reduce oscillations in the errors,
reduce errors and increase convergence rates.  SIAC has its basis in the work by
Bramble and Schatz \cite{BS} and Cockburn, Luskin, Shu and S\"{u}li
\cite{CLSS2}.  It has been extended to include boundary filtering and
nonuniform meshes \cite{ryan2003,SRV,Ryan1,Ryan2,SIACtriangles} as well as to
increase computational efficiency \cite{MRK1,SIACimplement}.  It has proven
effective for applications in visualization \cite{Steffanetal,OneStream} and it
is expected that this superconvergence extraction technique can act as a
natural filter for EFLES.

The traditional application of SIAC filtering takes the numerical
approximation, $u_h(x)$, and convolves it against a specially designed kernel,
\begin{equation}
u_h^\star(x)=\frac{1}{H}\int_{-\infty}^{\infty}\, K^{m+1,\ell}\left(\frac{y-x}{H}\right)u_h(y)\, dy,
\label{eq:siacdg}
\end{equation}
where $K^{m+1,\ell}$ is a linear combination of  $m$ B-splines of order $\ell$
and $H$ is a scaling typically related to a mesh quantity.

The symmetric form of the post-processing kernel can be written as
\begin{equation}
K^{m+1,\ell}(x)=\sum_{\gamma}\, c_\gamma^{m+1,\ell}\psi^{(\ell)}(x-\gamma),
\label{eq:kernel}
\end{equation}
where  $\psi^{(\ell)}(x)$ is obtained by convolving the characteristic function over the interval $(-\half, \half)$ with itself $\ell+1$ times and the coefficients $c_\gamma^{m+1,\ell}\in \mathbb{R}.$

The form of the SIAC kernel is similar to the form of the mixed polynomial \cite{Suarez2,Tornberg}.  Primary properties of the SIAC that make the SIAC kernel suitable for regularization include:
\begin{itemize}
  \item $\psi^{(\ell)}(x)$ can be expressed as a linear combination of Delta-functions, using the property $\frac{d^\alpha\psi^{(\ell)}(x)}{dx^\alpha}=\partial^\alpha_H\psi^{(\ell-\alpha)}(x),$ where $\partial_H$ represents a divided difference.
  \item The SIAC kernel can reproduce polynomials of degree $m.$  This is equivalent to the conditions
  \begin{itemize}
    \item $\int_\mathbb{R}\, K^{m+1,\ell}(\xi)\, d\xi=1.$
    \item $\int_\mathbb{R}\,  \xi^i K^{m+1,\ell}(\xi)\, d\xi=0$ for $i=1,\dots,m.$
  \end{itemize}
   These are similar to the conditions on mixed polynomials.  However, unlike the mixed polynomials, the SIAC kernel does not require the smoothness of the kernel even though it vanishes, at the endpoints of the compact support.
\end{itemize}

\subsection{Dirac-Delta Filtering}

\newcommand{\intdmk}{\int_{x-\varepsilon}^{x+\varepsilon} u(\tau) \delta_{\varepsilon}^{m,k}(x-\tau) d\tau}

In \cite{Suarez3}, an extension from \cite{Suarez2}, {\em singular source
terms} are expressed as weighted summations of Dirac-Delta kernels that are
regularized through approximation theory with convolution operators. The
regularization is obtained by convolution with the high-order compactly
supported Dirac-Delta kernel, (\ref{eq:delta_function}), whose overall accuracy
is controlled by the number of vanishing moments $m$, degree of smoothness $k$
and length of the support $\veps$.  In this work, the Dirac-Delta kernel is
used to regularize time dependent {\em discontinuous solutions}. Suppose the
solution is given by the variable $u(x), \, x \in[-1,1]$. Then the filtered
data, $\tilde{u}(x)$, follows from the convolution of $u(x)$ with the
Dirac-Delta kernel as
\begin{equation}
  \tilde{u}(x)=\int_{-1}^1 u(\tau) \, \dmk (x-\tau) d\tau.
\end{equation}
or, simply,
\begin{equation}
\tilde{u}(x)=\int_\xeps u(\tau) \, \dmk (x-\tau) d\tau,
\label{eq:convolution}
\end{equation}
since the Dirac-Delta kernel is zero outside its compact support $\xeps$.

To apply the filtering operation, we need to choose the number of vanishing
moments $m$, the number of continuous derivatives $k$ and the scaling parameter
$\varepsilon$. 

We make the following notes:
\begin{itemize}
  \item The number of vanishing moments and support width determine the accuracy of the Delta-kernel and therefore the error introduced by the filtering operation in smooth regions of the solution. In \cite{Suarez3} it was proven that the filtering error is $O(\varepsilon^{m+1})$, provided the scaling parameter is $\varepsilon=O(N^{(-k/m+k+2)})$. The requirement on the scaling parameter follows from the fact that the error in the quadrature rule used to evaluate the convolution integral has to be smaller than the $O(\varepsilon^{m+1})$ accuracy of the Dirac-Delta approximation. More vanishing moments reduces the filtering error, however it requires a wider support, leading to a wider regularization zone.

  \item The smoothness of the Dirac-Delta kernel is controlled by the number of continuous derivatives $k$ at the endpoints. In \cite{Suarez3} and \cite{Suarez2} it is shown that $k$ controls the smoothness of the transition between the regularized source and the solution and thus controls the order of convergence away from the source. When filtering the entire solution there is no such transition and thus the influence of $k$ is minor. This is confirmed by numerical experiments.
\end{itemize}

\begin{figure}[htbp]
\begin{center}
  \mbox{
   \makebox[0.40\textwidth][c]{$(m,k)=(3,8)$}
   \makebox[0.40\textwidth][c]{$(m,k)=(5,8)$}
  }
  \mbox{
    \includegraphics[width=0.40\linewidth]{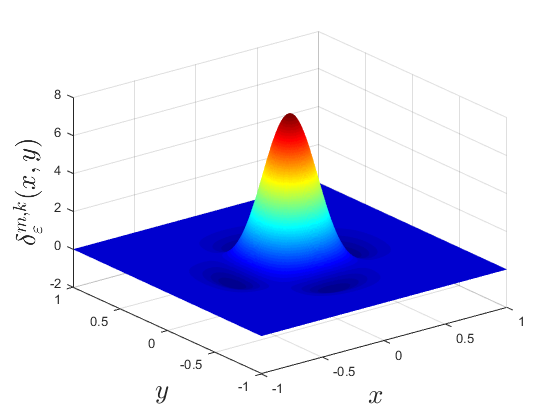}
    \includegraphics[width=0.40\linewidth]{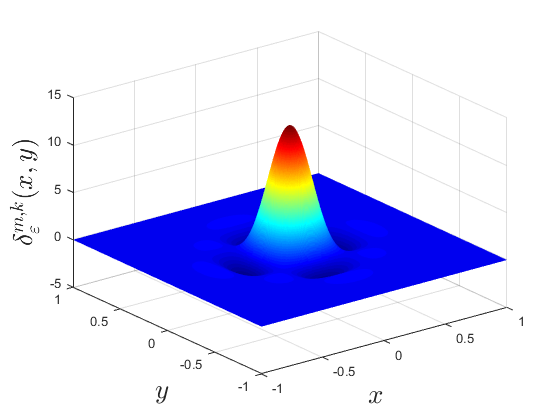}
   }
\end{center}
  \caption{(Color online) Tensor product of the one-dimensional Dirac-Delta function, for (Left) $(m,k)=(3,8)$ and (Right) $(m,k)=(5,8)$ with scaling parameter $\varepsilon=1$.}
\label{fig:kernel2D}
\end{figure}

Extension to two dimensions follows straightforwardly from the tensor product of the one-dimensional Delta-function:
\begin{equation}
\delta_{\varepsilon}^{m,k}(x,y)=\delta_{\varepsilon}^{m,k}(x) \otimes \delta_{\varepsilon}^{m,k}(y)
\label{eq:tensorproduct}
\end{equation}
Fig. \ref{fig:kernel2D} shows the two-dimensional equivalents of the one-dimensional kernels shown in Fig. \ref{fig:kernel1D}.
As in the one-dimensional case, the filtering is based on the convolution of the solution with the Delta-kernel,

\begin{equation}
\tilde{u}(x,y)=\int_{\Omega^\veps_x}\int_{\Omega^\veps_y} u(\tau,\eta) \delta_{\varepsilon}^{m,k}(x-\tau,y-\eta) d\tau d\eta.
\label{eq:convolution2D}
\end{equation}

\subsection{Implementation of the Filtering Operation}

Filtering of the interpolant of the solution (\ref{eq:interpolant}), leads to,
\begin{eqnarray}
\tilde{u}_N(x)=\int_\xeps \left[\sum\limits_{i=0}^N u(x_i) l_i(\tau)\right] \delta_{\varepsilon}^{m,k}(x-\tau)d\tau =\sum\limits_{i=0}^N u(x_i) S_i(x),
\label{eq:filtermatrix}
\end{eqnarray}
after interchanging the summation and integration, where the filtering function $S_i(x)$ is
\begin{equation}
  S_i(x) = \int_\xeps l_i(\tau) \delta_{\varepsilon}^{m,k}(x-\tau)d\tau.
\end{equation}
Hence, one has, at the collocation points,
\begin{equation}
  \overrightarrow{\tilde{u}} = \bS \vecu ,
\end{equation}
where the $(N+1)\times (N+1)$ filtering matrix $\bS$ has elements
\begin{equation}
  S_{i,j} = \int_\xepsj l_i(\tau) \delta_{\varepsilon}^{m,k}(x_j-\tau)d\tau.
	\label{filtermatrix_element}
\end{equation}

In two dimensions, the filtering operation can be written compactly as
\begin{equation}
  \tilde{\bf U} = \bS_x \bf U \bS_y^\ssT,
\label{filtering2D}
\end{equation}
where $\bS_x$ and $\bS_y$ are the one-dimensional filtering matrix in $x-$ and
$y-$ direction respectively and the superscript $\ssT$ denotes transpose.

\begin{remark}
Near the boundaries, the compact support of the high-order Dirac-Delta function, $\delta_{\varepsilon}^{m,k}$, extends out of the domain and hence, the data can not be filtered in that case.  This means that $\tilde{u}(x)=u(x)$ for $|x| > 1-\veps$. The filtering matrix $S$ can be precomputed and stored for later use as long as the filter parameters remain unchanged.
\end{remark}

\subsection{Clenshaw-Curtis quadrature}

The one remaining important issue that needs to be addressed is how to evaluate
the integrals in Eq. \ref{eq:filtermatrix}. Since the high-order Dirac-Delta function is a polynomial of degree $M=m+2(k+1)$ and the Lagrange interpolation
polynomial is of degree $N-1$, the integrand is a polynomial of degree $M+N-1$
and can be evaluated analytically. Since this can become time-consuming for
large $N$ an appropriate quadrature rule is preferred. For this purpose
Clenshaw-Curtis quadrature will be used, that is, \begin{equation}
  \int_\xepsi l_n(\tau) \delta_{\varepsilon}^{m,k}(x_i-\tau)d\tau \approx \sum\limits_{q=0}^Q w_q l_n(x_q) \delta_{\varepsilon}^{m,k}(x_i-x_q).
\label{eq:clenshaw}
\end{equation}
where $Q$ is the number of Chebyshev Gauss-Lobatto quadrature nodes used in the compact support domain $\xepsi$, that is,
\begin{equation}
   x_q=x_i-\varepsilon \cos \left( \frac{\pi q}{Q} \right), \quad q=0,\ldots,Q,
\end{equation}
and $w_q$ are the corresponding weights. Clenshaw Curtis quadrature exactly
evaluates polynomials of degree $Q-1$. Hence, if one takes $Q=M+N$, the
integrals are evaluated exactly. Also, the weights $w_q$ can be precomputed
using fast Fourier transform (FFT).


\subsection{Shock Capturing with Dirac-Delta Filtering}

The theoretical requirement for the scaling parameter,
$\varepsilon=\mathcal{O}(N^{(-k/m+k+2)})$, to assure
$\mathcal{O}(\varepsilon^{m+1})$ accuracy in the convolution operation is based
on the requirement that the error in the quadrature rule has to be smaller than
the error in the Dirac-Delta approximation \cite{Suarez3}.   Since the
convolution integrals in this work are solved exactly using Clenshaw-Curtis
quadrature, this criterium can be relaxed. For shock capturing we choose the
scaling parameter to be proportional to the grid spacing and to guarantee that
a minimal of at least two neighboring collocation points will be located inside
the compact support width $\xeps_i$ for any point $x=x_i$.  In this case, the
scaling parameter will be expressed in terms of the number of points $N_d$ the
kernel spans at the center of the domain $x=x_{N/2}$, assuming $N$ is even,
that is, \begin{equation}
  \varepsilon=\sin \left(\pi N_d/(2N) \right).
\end{equation}
The value of $N_d$ is determined empirically and chosen to ensure stable
converging results.  It depends on the number of vanishing moments $m$ of the
Dirac-Delta kernel and whether the filter is applied in a linear or non-linear
equation.  For implementation in the linear advection equation, only the
initial condition is filtered while for a nonlinear PDE's, such as the Burgers'
equation and the Euler equations, the solution is filtered at every Runge Kutta
time step. This is because of the formation of a finite space-time singularity
by the nonlinearity of the equations.
%
%
However, filtering can lead to smearing of the discontinuity and the summation
of filter errors in smooth regions.  In order to limit these less desirable
effects, we use a smaller scaling parameter, $\varepsilon$, as compared to the
linear advection equation.

\begin{remark}

Near the boundary the filter cannot be convoluted due to the symmetry of the convolution kernel.
	In this paper those few points are reset to be the exact solution for the nonlinear PDEs in order
to isolate and study solely the effects of the Dirac-Delta kernels.
Furthermore, to avoid the effect of variable time step $\dt$, we fixed the
stable time step of $\dt=1 \times 10^{-5}$ in all the simulations performed
below in order to be able to compare results.  \end{remark}



\section{Numerical Tests} \label{chap:Results}

\subsection{Advection equation} \label{Advection equation}

We shall consider the simple one-dimensional linear advection equation with a discontinuity as initial condition on the domain $-1<x<1$,
\begin{eqnarray}
\frac{\partial u}{\partial t} + \frac{\partial u}{\partial x} &=& 0, \\
u(x,0)&=&\left\{
        \begin{array}{ll}
           \sin(\pi x) -0.5 \qquad x \leq -0.25 \\
           \sin(\pi x) +0.5 \qquad x > -0.25
         \end{array} \right. ,  \nonumber \\
u(-1,t)&=&\sin(\pi(-1-t)) -0.5. \nonumber
\label{eq:advection2}
\end{eqnarray}
The initial condition is first filtered by the filtering matrix $\bS$.
For the results below we used a Dirac-Delta kernel with $(m,k)=(3,8)$ and $N_d=13$.
\begin{figure}[htbp]
\begin{center}
  \mbox{
        \includegraphics[width=0.40\linewidth]{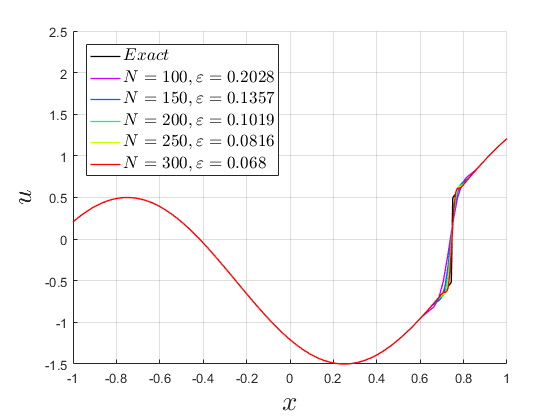}
        \includegraphics[width=0.40\linewidth]{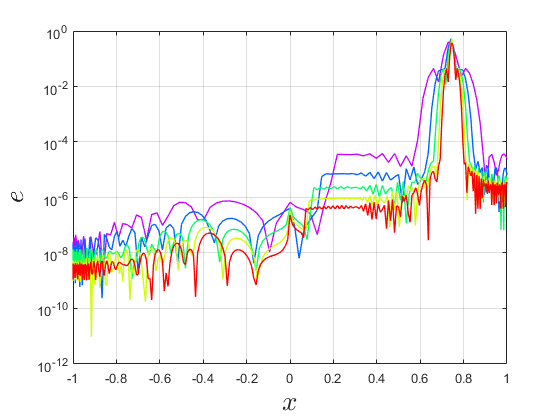}
       }
\end{center}
    \caption{Linear advection equation: Analytical and spectral solution (Left), and pointwise error (Right) at time $t=1$, for the five different grids with $(m,k)=(3,8)$ and $N_d=13$.}
\label{fig:advection0}
\end{figure}

\begin{figure}[htbp]
\begin{center}
  \mbox{
        \includegraphics[width=0.40\linewidth]{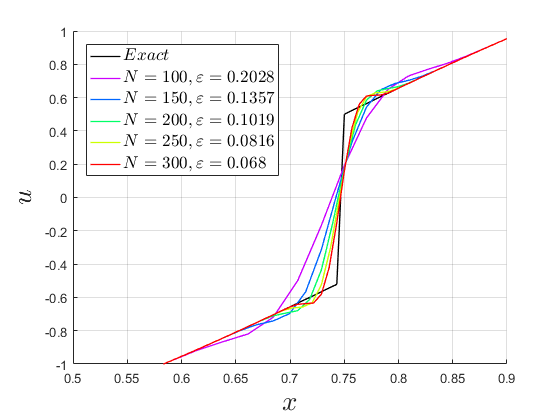}
        \includegraphics[width=0.40\linewidth]{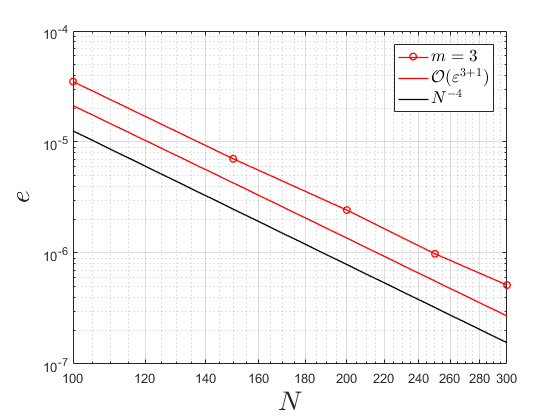}
       }
\end{center}
    \caption{Linear advection equation: Detail view of regularization zone (Left) and error convergence at $x=0.28$ compared to the theoretical value $O(\varepsilon^{m+1})$ (Right) for the five different grids with $(m,k)=(3,8)$ and $N_d=13$.}
\label{fig:advection1}
\end{figure}

In Fig. \ref{fig:advection0}, the point-wise error clearly shows the advected filter error for $x>0$ caused by filtering the initial condition.  In Fig. \ref{fig:advection1}, the rate of convergence outside the regularization zone is examined at $x=0.28$ and follows the theoretical value of $O(\varepsilon^{m+1})$.

\subsection{Burger's equation} \label{burgers_equation}

Next, we consider the Burgers' equation,
\begin{eqnarray}
\frac{\partial u}{\partial t} + \frac{1}{2}\frac{\partial u^2}{\partial x} &=& 0, \\
u(x,0)&=& -\sin(\pi x), \\
u(\pm 1,t)&=&0 \nonumber
\label{eq:burgers2}
\end{eqnarray}
In this case, the solution is filtered at the end of every Runge-Kutta time step.
For the results below we used a kernel with $(m,k)=(3,8)$ and $N_d=2.5$.
\begin{figure}[htbp]
\begin{center}
  \mbox{
        \includegraphics[width=0.40\linewidth]{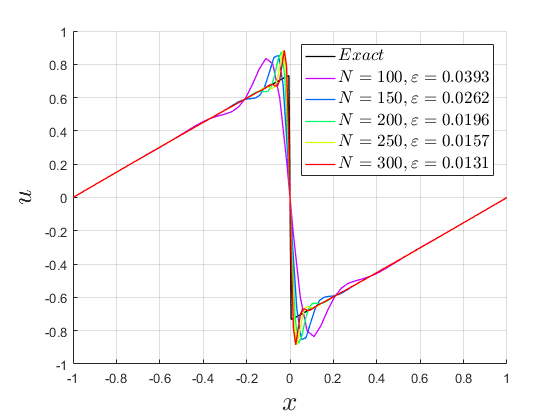}
        \includegraphics[width=0.40\linewidth]{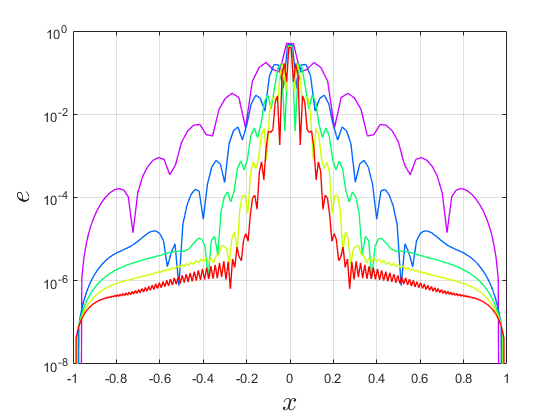}
       }
\end{center}
    \caption{Burger's equation: Analytical and spectral solution (Left), and pointwise error (Right) at time $t=1$, for the five different grids with $(m,k)=(3,8)$ and $N_d=2.5$.}
\label{burgers}
\end{figure}

The solution in  Figure \ref{burgers} shows that the
filtering effectively captures the shock and gives a stable solution. A clearly
visible overshoot is however introduced. Because of the non-linear character of the
solution, increasing the number of vanishing moments does improve
the result. 

The error (Fig. \ref{burgers}) that is introduced by regularization at the shock
location is smeared over the domain by the  multiple application of the filter
at each time step.
For the values of $N_d$ used here, the smearing is is small.
The smearing effect increases with increasing support width. 
We not that extending the simulation time until $t=5$ does
not change the error much as compared to error at $t=1$.
Only the shock is slightly more dissipated.

\subsection{Euler equations} \label{euler_equations}

The 2D unsteady Euler Equations in conservative form are given as

\begin{equation}
\frac{\partial U}{\partial t} + \frac{\partial F_x}{\partial x} + \frac{\partial F_y}{\partial y}= 0,
\label{euler_1}
\end{equation}
where the conserved quantities $U$ and the flux vectors $F_x$ and $F_y$ are given by
\begin{equation}
U=
\begin{pmatrix}
\rho \\ \rho u \\ \rho v \\ \rho E
\end{pmatrix}
,
F_x=
\begin{pmatrix}
\rho u \\ \rho u^2 + p \\ \rho u v \\ (\rho E + p) u
\end{pmatrix}
,
F_y=
\begin{pmatrix}
\rho u \\ \rho u v\\ \rho v^2 + p \\ (\rho E + p) v
\end{pmatrix}.
\label{euler_2}
\end{equation}
The system is closed by assuming a calorically perfect gas, which relates the pressure to the density, velocity and energy. For the two-dimensional case this gives
\begin{equation}
p=(\gamma -1)\left(\rho E - \frac{1}{2}\rho (u^2 +v^2) \right),
\label{euler_3}
\end{equation}
in which $\gamma$ is the specific heat ratio, which is $1.4$ for air. 

\subsubsection{Sod's shock tube problem}

Sod's shock tube problem is governed by the 1D Euler equations. The initial conditions for Sod's shock tube problem are
\begin{table}[H]
\centering
\label{my-label}
\begin{tabular}{ll}
$\rho_L=1.0$      & $\rho_R=0.125,$ \\
$p_L=1.0$         & $p_R=0.1,$                 \\
$u_L=0.0$         & $u_R=0.0.$
\end{tabular}
\end{table}
where the left state is in the domain with $x < 0$ and the right state is in
$x \geq 0$
The conserved quantities are filtered every time-step using the filtering matrix. For the results below we used a kernel with $(m,k)=(3,8)$ and $N_d=2.5$.
\begin{figure}[htbp]
\begin{center}
  \mbox{
        \includegraphics[width=0.40\linewidth]{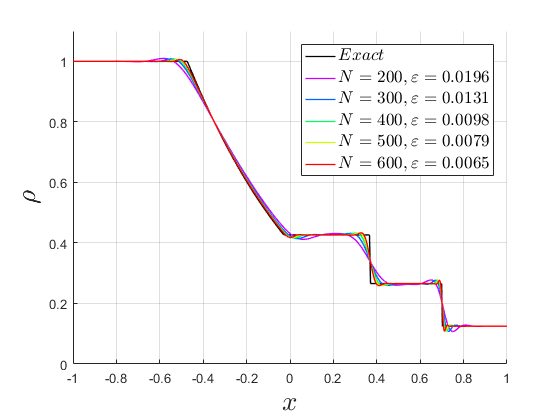}
        \includegraphics[width=0.40\linewidth]{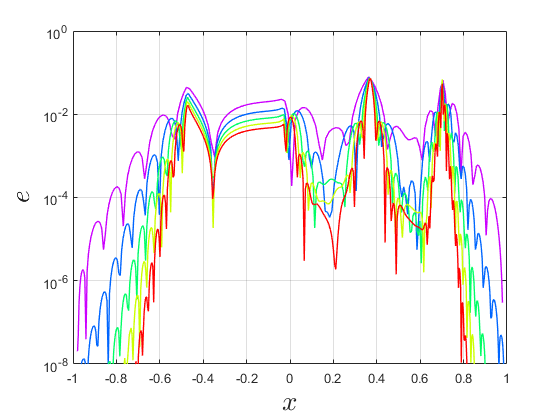}
       }
\end{center}
    \caption{Analytical solution, spectral solution (Left) and pointwise error (Right) of the density at $t=0.4$, for the five different grids for $(m,k)=(3,8)$ and $N_d=2.5$.}
\label{fig:sod1D}
\end{figure}
The solutions show that both the shock and the contact discontinuity are effectively captured using the filter-operation. 
A low error is observed outside the regularization zone.

\subsubsection{Shu-Osher problem}

Shu-Osher's Problem is also governed by the 1D Euler equations with the initial conditions of the primitive variables in the left state $(\rho_L, P_L, u_L) = (27/7, 31/3, 4\sqrt{35/9}$ and the right state $(\rho_R, P_R, u_R) = (1+0.2\sin(25x), 1, 0)$.
Again, the conserved quantities are filtered every time-step using the
filtering matrix. For the results below we used a kernel with $(m,k)=(5,8)$ and
$N_d=5.5$ in order to limit the error introduced in the entropy-wave. A $5$th
order WENO-solution on $10000$ points is used to serve as a reference 'exact'
solution.

\begin{center}
\large{\textbf{Results, $\delta_{\varepsilon}^{5,8}$, $N_d=6.5$}}
\end{center}

\begin{figure}[H]
\begin{center}
  \mbox{
        \includegraphics[width=0.40\linewidth]{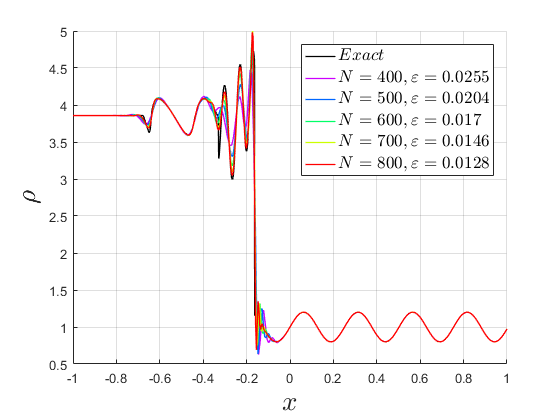}
        \includegraphics[width=0.40\linewidth]{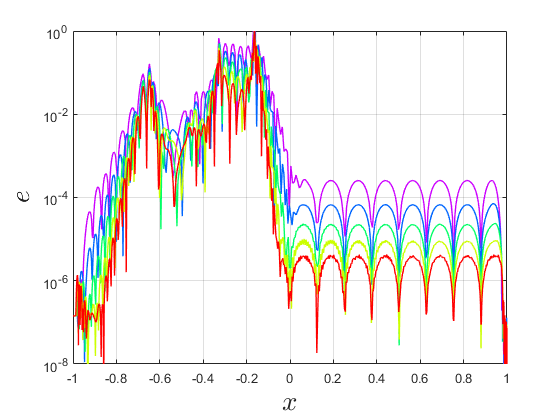}
       }
\end{center}
    \caption{Analytical solution, spectral solution (Left) and pointwise error (Right) of the density at $t=0.18$, for the five different grids for $(m,k)=(5,8)$ and $N_d=6.5$.}
\end{figure}

\begin{figure}[H]
\begin{center}
  \mbox{
        \includegraphics[width=0.40\linewidth]{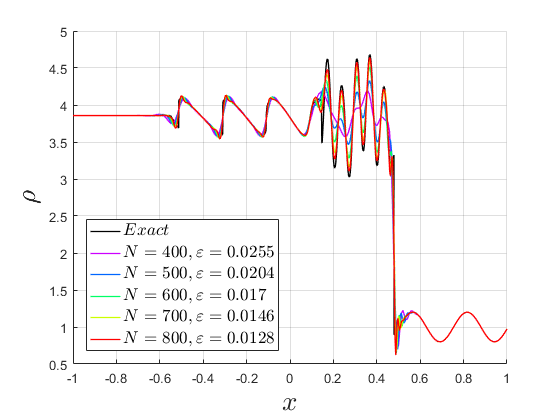}
        \includegraphics[width=0.40\linewidth]{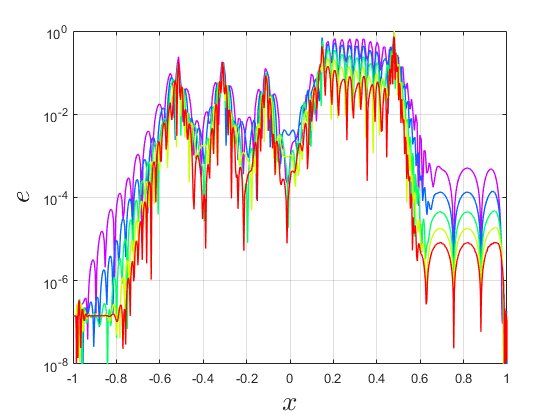}
       }
\end{center}
    \caption{Analytical solution, spectral solution (Left) and pointwise error (Right) of the density at $t=0.36$, for the five different grids for $(m,k)=(5,8)$ and $N_d=6.5$.}
\label{fig:shuosher0}
\end{figure}

\begin{figure}[H]
\begin{center}
  \mbox{
        \includegraphics[width=0.40\linewidth]{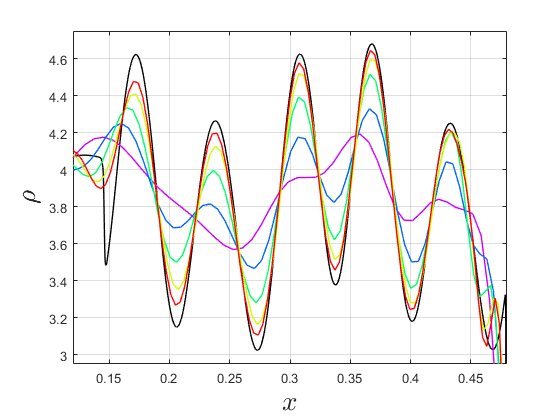}
        \includegraphics[width=0.40\linewidth]{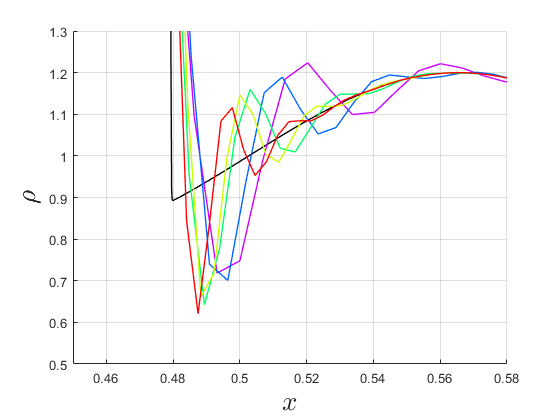}
       }
\end{center}
    \caption{Detail view of the entropy wave (Left) and the regularization zone (Right) at $t=0.36$, for the five different grids for $(m,k)=(5,8)$ and $N_d=6.5$.}
\label{fig:shuosher1}
\end{figure}




The solution shows that the filter effectively captures the shock while the kernel with $m=5$ vanishing moments ensures high resolution behind the shock. Fig. \ref{fig:shuosher1} provides a detail view of the entropy wave (Left) and the regularization zone (Right).

\subsection{Explosion problem}

As a 2D test case, the explosion problem is considered. This problem is
governed by the 2D Euler equations, (\ref{euler_1}), (\ref{euler_2}) and
(\ref{euler_3}). The equations are solved on a $2.0 \times 2.0$ square domain
in the $x-y$ plane. The initial condition consists of the region inside of a
circle with radius $R=0.4$ centered at $(0,0)$ and the region outside of the
circle. The flow variables are constant in each of these regions and are
separated by a circular discontinuity at time $t=0$. The two constant states
are chosen as: \begin{eqnarray}
\rho_{in}&=&1.0 \quad \rho_{out}=0.125, \\
p_{in}&=&1.0 \quad p_{out}=0.1, \\
u_{in}&=&0.0 \quad u_{out}=0.0, \\
v_{in}&=&0.0 \quad v_{out}=0.0,
\end{eqnarray}
where the subscripts $in$ and $out$ denote values inside and outside the circle, respectively. The approach described by Eleuterio \cite{Eleuterio} is used to serve as an exact solution. For the results below we used a kernel with $(m,k)=(3,8)$ and $N_d=2.5$.
\begin{figure}[htbp]
\begin{center}
  \mbox{
        \includegraphics[width=0.40\linewidth]{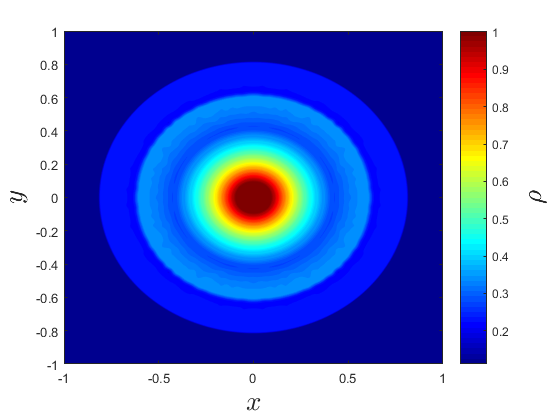}
       }
\end{center}
\caption{Contour plot of the density at $t=0.25$ for $(m,k)=(3,8)$ and $N_d=2.5$.}
\label{fig:explosion0}
\end{figure}

\begin{figure}[htbp]
\begin{center}
  \mbox{
        \includegraphics[width=0.40\linewidth]{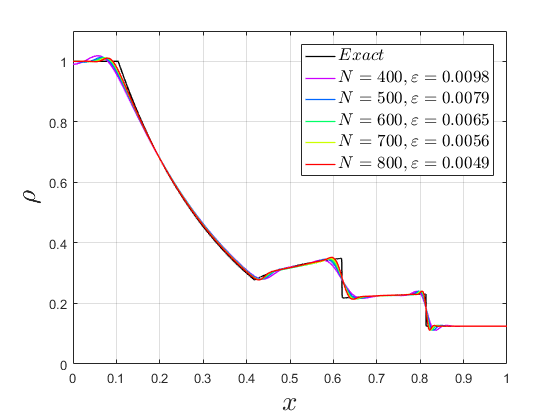}
        \includegraphics[width=0.40\linewidth]{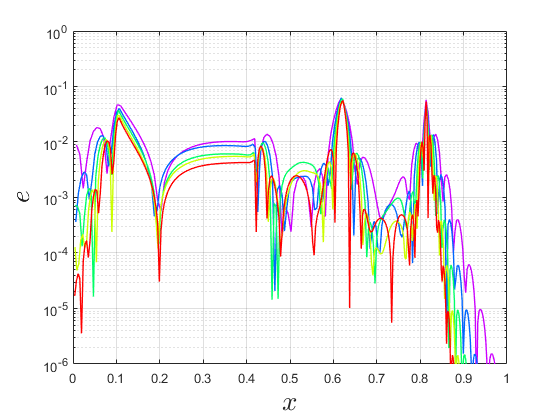}
       }
\end{center}
    \caption{Numerical solution (Left) and error (Right) along the radial line at $t=0.25$, for the five different grids for $(m,k)=(3,8)$ and $N_d=2.5$.}
\label{fig:explosion1}
\end{figure}
The contour plot of the density shows that the symmetry of the flow is captured
on the 2D Cartesian grid. Comparisons along any radial line give identical
results.

\section{Conclusions and future work} \label{chap:Conclusions}

In this paper, the use of high-order Dirac-Delta function based filters for the
regularization of shocks and discontinuities in combination with a global
spectral method is investigated. We have shown that these filters are able to
effectively capture shocks while maintaining high resolution in smooth areas of
the solution.

As a first next step for future work, the filtering operation has to be
extended for the use in discontinuous spectral element methods. Furthermore the
development of a way to apply the filtering near the boundaries is desired. A
possible way could be the development of special one-sided polynomial kernels,
analog to the one-sided SIAC kernels for boundary filtering \cite{OneStream},
\cite{ryan2003}.

In the present work, the developed filters have been applied globally, i.e. the
whole domain is affected by the filter, thus causing (small) errors in smooth
areas and smearing of the discontinuity. It would be beneficial if this effect
could be limited. This could be accomplished if only the area of the
discontinuity is filtered. This however requires some means to identify shocks
and discontinuities in the solution. One possibility is the use of
multi-resolution analysis (MR) \cite{Harten1}, \cite{Harten2}. MR-analysis is
successfully used in hybrid spectral-WENO methods to switch between
high-resolution WENO methods for shocked regions and computationally efficient
spectral methods in smooth regions \cite{hybridWENO}.
In a follow on paper, we intend to report on the 
use the filter-error itself to identify shocks
and discontinuities. 
The filter-error \begin{equation}
  e=|\boldsymbol{u}- \bS \boldsymbol{u}|,
\end{equation}
describes the shape of the solution in the sense that it is small in smooth areas and large in shocked areas. If a function $L(e)$ can be derived that is able to determine which areas should be filtered and which not, the filtering operation could be improved and implemented according to
\begin{equation}
  \tilde{\boldsymbol{u}}=L(e)\boldsymbol{u} + (1-L(e))\bS \boldsymbol{u}.
\end{equation}
Thus for areas in which $L(e)=1$ the solution is not filtered whereas if $L(e)=0$ the solution is fully filtered. This would then effectively limit the smearing and the error introduced in smooth areas.


\bibliographystyle{unsrt}
\bibliography{biblio}

\begin{thebibliography}{10}

\bibitem{Shu}
C.W. Shu and S.~Osher.
\newblock Efficient {I}mplementation of {E}ssentially {N}on-oscillatory
  {S}hock-{C}apturing {S}chemes.
\newblock {\em Journal of Computational Physics}, 77:439--471, 1988.

\bibitem{Cockburn}
B.~Cockburn.
\newblock Devising discontinuous galerkin methods for non-linear hyperbolic
  conservation laws.
\newblock {\em J. Comp. Apll. Math.}, 128:187--204, 2001.

\bibitem{Cockburn2}
B.~Cockburn and C.W. Shu.
\newblock Runge-kutta discontinuous galerkin methods for convection-dominated
  problems.
\newblock {\em Journal of scientific computing}, 16:173--261, 2001.

\bibitem{Chaudhuri}
A.~Chaudhuri, G.B. Jacobs, W.S. Don, H.~Abassi, and F.~Mashayek.
\newblock Explicit discontinuous spectral element method with entropy
  generation based artificial viscosity for shocked viscous flows.
\newblock {\em J. Comp. Phys.}, 332:99--117, 2017.

\bibitem{Hughes}
T.J. Hughes, L.~Franca, and M.~Mallet.
\newblock A new finite element formulation for computational fluid dynamics: I.
  symmetric forms of the compressible euler and navier-stokes equations and the
  second law of thermodynamics.
\newblock {\em Computer Methods in Applied Mechanics and Engineering},
  54:223--234, 1986.

\bibitem{Guermond}
J.L. Guermond, R.~Pasquetti, and B.~Popov.
\newblock Entropy viscosity method for nonlinear conservation laws.
\newblock {\em Journal of Computational Physics}, 230:4248--4267, 2011.

\bibitem{don}
W.S. Don.
\newblock Numerical {S}tudy of {P}seudospectral {M}ethods in {S}hock {W}ave
  {A}pplications.
\newblock {\em J. Comput. Phys.}, 110:103--111, 1994.

\bibitem{vandeven}
H.~Vandeven.
\newblock Family of spectral filters for discontinuous problems.
\newblock {\em Journal of Sientific Computing}, 6(2):159--192, 1991.

\bibitem{BS}
J.H. Bramble and A.H. Schatz.
\newblock Higher {O}rder {L}ocal {A}ccuracy by {A}veraging in the {F}inite
  {E}lement {M}ethod.
\newblock {\em Mathematics of Computation}, 31:94--111, 1977.

\bibitem{CLSS2}
B.~Cockburn, M.~Luskin, C.W. Shu, and E.~Süli.
\newblock Enhanced {A}ccuracy by {P}ost-processing for {F}inite {E}lement
  {M}ethods for {H}yperbolic {E}quations.
\newblock {\em Mathematics of Computation}, 72:577--606, 2003.

\bibitem{MockLax}
M.S. Mock and P.D. Lax.
\newblock The computation of discontinuous solutions of linear hyperbolic
  equations.
\newblock {\em Comm. Pure Appl. Math.}, 18:423--430, 1978.

\bibitem{Suarez2}
J.P. Suarez, G.B. Jacobs, and W.S. Don.
\newblock A high-order {D}irac-delta regularization with optimal scaling in the
  spectral solution of one-dimensional singular hyperbolic conservation laws.
\newblock {\em SIAM J. Sci. Comp}, 36:1831--1849, 2014.

\bibitem{Tornberg}
A.K. Tornberg.
\newblock Multi-dimensional quadrature of singular and discontinuous functions.
\newblock {\em BIT}, 42:644--669, 2002.

\bibitem{Suarez3}
J.P. Suarez and G.B. Jacobs.
\newblock Regularization of singularities in the weighted summation of
  {D}irac-delta functions for the spectral solution of hyperbolic conservation
  laws.
\newblock {\em J. Sci. Comp.}, \textit{preprint at:
  https://arxiv.org/abs/1611.05510}.

\bibitem{Hesthaven}
Sigal~Gottlieb Jan~Hesthaven and David Gottlieb.
\newblock Spectral {M}ethods for {T}ime-{D}ependent {P}roblems.
\newblock {\em Cambridge University Press, Cambridge, United Kingdom}.

\bibitem{RSA}
J.K. Ryan, C.W. Shu, and H.L. Atkins.
\newblock Extension of a post-processing technique for discontinuous {G}alerkin
  methods for hyperbolic equations with application to an aeroacoustic problem.
\newblock {\em SIAM Journal on Scientific Computing}, 26:821--843, 2004.

\bibitem{kingmirzaee}
J.~King, H.~Mirzaee, J.K. Ryan, and R.M. Kirby.
\newblock Smoothness-{I}ncreasing {A}ccuracy-{C}onserving ({SIAC}) filtering
  for discontinuous {G}alerkin solutions: {I}mproved errors versus higher-order
  accuracy.
\newblock {\em Journal of Scientific Computing}, 53:129--149, 2012.

\bibitem{CKRS}
S.~Curtis, R.M. Kirby, J.K. Ryan, and C.W. Shu.
\newblock Post-processing for the discontinuous {G}alerkin method over
  non-uniform meshes.
\newblock {\em SIAM Journal on Scientific Computing}, 30:272--289, 2007.

\bibitem{JIRyan}
L.~Ji, P.~van Slingerland, J.K. Ryan, and C.~Vuik.
\newblock Superconvergent error estimates for a position-dependent
  {S}moothness-{I}ncreasing {A}ccuracy-{C}onserving filter for {DG} solutions.
\newblock {\em Mathematics of Computation}, 83:2239--2262, 2014.

\bibitem{ryan2003}
J.K. Ryan and C.W. Shu.
\newblock On a one-sided post-processing technique for the discontinuous
  {G}alerkin methods.
\newblock {\em Methods Appl. Anal.}, 10:295--307, 2003.

\bibitem{SRV}
P.~van Slingerland, J.K. Ryan, and C.~Vuik.
\newblock Position-dependent {S}moothness-{I}ncreasing {A}ccuracy-{C}onserving
  ({SIAC}) filtering for accuracy for improving discontinuous {G}alerkin
  solutions.
\newblock {\em SIAM J. Sci. Comp.}, 33:802--825, 2011.

\bibitem{Ryan1}
X.~Li, R.M. Kirby, and C.~Vuik.
\newblock Computationally efficient position-dependent
  {S}moothness-{I}ncreasing {A}ccuracy-{C}onserving ({SIAC}) filtering: {T}he
  uniform mesh case.
\newblock {\em Preprint}, 2013.

\bibitem{Ryan2}
X.~Li, R.M. Kirby, and C.~Vuik.
\newblock Computationally efficient position-dependent
  {S}moothness-{I}ncreasing {A}ccuracy-{C}onserving ({SIAC}) filtering: {T}he
  nonuniform mesh case.
\newblock {\em Preprint}, 2014.

\bibitem{SIACtriangles}
H.~Mirzaee, L.~Ji, J.K. Ryan, and R.M. Kirby.
\newblock Smoothness-{I}ncreasing {A}ccuracy-{C}onserving ({SIAC})
  post-processing for discontinuous {G}alerkin solutions over structured
  triangular meshes.
\newblock {\em SIAM Journal on Numerical Analysis}, 49:1899--1920, 2011.

\bibitem{MRK1}
H.~Mirzaee, J.K. Ryan, and R.M. Kirby.
\newblock Quantification of errors introduced in the numerical approximation
  and implementation of {S}moothness-{I}ncreasing {A}ccuracy-{C}onserving
  ({SIAC}) filtering of discontinuous {G}alerkin ({DG}) fields.
\newblock {\em Journal of Scientific Computing}, 45:447--470, 2010.

\bibitem{SIACimplement}
H.~Mirzaee, J.K. Ryan, and R.M. Kirby.
\newblock Efficient implementation of {S}moothness-{I}ncreasing
  {A}ccuracy-{C}onserving ({SIAC}) filters for discontinuous {G}alerkin
  solutions.
\newblock {\em Journal of Scientific Computing}, 52:85--112, 2012.

\bibitem{Steffanetal}
M.~Steffan, S.~Curtis, R.M. Kirby, and J.K. Ryan.
\newblock Investigation of {S}moothness-{I}ncreasing {A}ccuracy-{C}onserving
  {F}ilters for {I}mproving {S}treamline {I}ntegration through {D}iscontinous
  {F}ields.
\newblock {\em IEEE-TVCG}, 14:680--692, 2008.

\bibitem{OneStream}
D.~Walfisch, J.K. Ryan, R.M. Kirby, and R.~Haimes.
\newblock One-sided {S}moothness-{I}ncreasing {A}ccuracy-{C}onserving filtering
  for enhanced streamline integration through discontinuous fields.
\newblock {\em Journal of Scientific Computing}, 38:164--184, 2009.

\bibitem{Eleuterio}
Eleuterio~F. Toro.
\newblock Riemann {S}olvers and {N}umerical {M}ethods for {F}luid {D}ynamics.
\newblock {\em Springer-Verlag, Berlin}, 2009.

\bibitem{Harten1}
A.~Harten.
\newblock High resolution schemes for hyperbolic conservation laws.
\newblock {\em Comput. Phys.}, 49:357--393, 1983.

\bibitem{Harten2}
A.~Harten.
\newblock Adaptive multiresolution schemes for shock computations.
\newblock {\em Comput. Phys.}, 115:319--338, 1994.

\bibitem{hybridWENO}
B.~Costa and W.S. Don.
\newblock Multi-domain hybrid spectral-{WENO} methods for hyperbolic
  conservation laws.
\newblock {\em Journal of Computational Physics}, 224:970--991, 2007.

\end{thebibliography}


\end{document}